\documentclass{amsart}

\newtheorem{theorem}{Theorem}[section]
\newtheorem{lemma}[theorem]{Lemma}
\newtheorem{corollary}[theorem]{Corollary}
\newtheorem{proposition}[theorem]{Proposition}

\theoremstyle{definition}

\theoremstyle{remark}
\newtheorem{remark}[theorem]{Remark}

\numberwithin{equation}{section}

\def\operator#1#2{\def#1{\mathop{\kern0pt\fam0#2}\nolimits}} 
\operator\ini{in}
\operator\Gin{Gin} 
\operator\GL{GL}
\operator\Support{Support}
\operator\mm{{\mathbf m}}
\operator\revlex{RL}
\operator\depth{depth}
\operator\Lex{Lex}
\operator\chara{char}

\begin{document}

\title{Reduction numbers and initial ideals}

\author{Aldo Conca}
\address{Dipartimento di
Matematica, Universit\'a di Genova, Via Dodecaneso 35, I-16146
Genova, Italia }
\email{ conca@dima.unige.it}




\begin{abstract}
The reduction number $r(A)$ of a standard graded algebra $A$ is the least integer $k$ such that
there exists a minimal reduction $J$ of the homogeneous maximal ideal $\mm$ of $A$
such that $J\mm^k=\mm^{k+1}$.
Vasconcelos conjectured  that  $r(R/I)\leq r(R/\ini(I))$  where
$\ini(I)$ is the initial ideal of an ideal $I$  in a polynomial ring $R$ with respect to a
term order.  The goal of this note is to prove the conjecture.  
\end{abstract}

\maketitle

\section{Reduction numbers and initial ideals}
Let $K$ be an infinite  field and let $A=\oplus_{i\in \bf N} A_i$  be a
homogeneous $K$-algebra, that is, an algebra of the form $R/I$   where
$R=K[x_1,\dots,x_n]$ is a polynomial ring and $I$ is a homogeneous ideal. 
The reduction number $r(A)$ of $A$ is the least integer $k$ such that there
exists a minimal reduction $J$ of the homogeneous maximal ideal $\mm$ of $A$
such that $J\mm^k=\mm^{k+1}$. It is not difficult to see that  $r(A)$ is
the  largest integer $k$ such that the Hilbert function of $A/J$ at $k$ does
not vanishes; here $J$ is the ideal of $A$ generated by   $d=\dim A$ generic
linear forms. 

Vasconcelos conjectured \cite[Conjecture 7.2]{V} that 
\[
r(R/I)\leq r(R/\ini_\tau(I))
\]
  where $\ini_\tau(I)$ is the initial ideal
of $I$ with respect to a term order $\tau$. The conjecture has been proved
by Bresinsky and Hoa \cite{BH} for the generic initial ideal $\Gin_\tau(I)$, or,
more generally, when $\ini_\tau(I)$ is Borel-fixed. Trung \cite{T} showed that  
$r(R/I)=r(R/\Gin_{\revlex}(I))$ where the $\Gin_{\revlex}(I)$ is the generic
initial ideal of $I$    with respect to the degree reverse lexicographic 
order $\revlex$ (revlex for short).  

The goal of this note is to prove the conjecture in general.  After this
paper was written we were informed  that Trung  \cite{T1} has independently
solved the conjecture in general by a completely different method.   What 
we prove is the following  generalization of Vasconcelos' conjecture: 
 
\begin{theorem}
\label{thm1}
   Let $p$ be an  integer, $0\leq p\leq n$,  and
let $\ini_\tau(I)$ be the initial ideal of $I$ with respect to a term order
$\tau$.  Let
$J$ be an ideal generated by $p$ generic linear forms. Then the Hilbert
function of $R/I+J$ is $\leq $  that  of  $R/\ini_\tau(I)+J$, that is 
\[
\dim_K \left [R/I+J \right ]_j \leq \dim_K \left [R/\ini_\tau(I)+J \right
]_j
\]  
for all $j\in {\bf N}.$  
\end{theorem}

Taking $p=\dim R/I$ and $j=r(R/\ini_\tau(I))+1$ one obtains
$r(R/I)<j$ which implies Vasconcelos' conjecture.   To prove Theorem \ref{thm1} we need some
preparation.

\begin{lemma}
\label{lm2}
  Let $p$ be an  integer, $0\leq p\leq n$,  and
let $J$ be an ideal generated by
$p$ generic linear forms. Then the Hilbert function of $R/I+J$ is equal to
the Hilbert function of $R/\Gin_{\revlex}(I)+H$ where  
$\Gin_{\revlex}(I)$  is the revlex Gin   of $I$ and
$H=(x_{n-p+1},x_{n-p+2},\dots,x_n)$.
\end{lemma}

\begin{proof} Set $J=(y_1,\dots,y_p)$. We take a matrix $g\in \GL_n(K)$
such that the induced $K$-algebra graded  homomorphism 
$g:R\to R$ maps $y_i$ to $x_{n-p+i}$ for $i=1,\dots,p$. It follows that the
Hilbert function of $R/I+J$  equals that of $R/g(I)+H$. Taking  initial
ideals does not change the Hilbert function and by  the known properties of
the revlex order    one has     
$\ini_{\revlex}(g(I)+H)=\ini_{\revlex}(g(I))+H$. But since the
$y_i$ are generic, $g$ is   generic as well. Then
$\ini_{\revlex}(g(I))=\Gin_{\revlex}(I)$ and we are done.  
\end{proof}

We would like now to  compare $\Gin_{\revlex}(I)$ with
$\Gin_{\revlex}(\ini_\tau(I))$. To this end, let us introduce a piece of
notation.  

Let $V\subset R_i$ be a subspace of forms of degree $i$ and dimension $d$. 
Then $\wedge^d V$ is a  subspace of dimension $1$ of  $\wedge^d R_i$. We
identify in the following  $\wedge^d V$ with any  non-zero element contained
in it.  Fix a term order $<_\sigma$ in $R$. An exterior monomial is an
element of the form $m_1\wedge\dots \wedge m_d$ where the $m_j$ are distinct
monomials of $R_i$.  An exterior monomial  $m_1\wedge\dots \wedge m_d$ is  
$\sigma$-standard    if  $m_1>_\sigma \dots>_\sigma  m_d$.  Note that
$\wedge^d R_i$ has a basis consisting of the $\sigma$-standard exterior
monomials.  We order the $\sigma$-standard  exterior monomials 
lexicographically:
\[
m_1\wedge\dots \wedge m_d>_\sigma  n_1\wedge\dots \wedge n_d
\]
 if
$m_i>_\sigma n_i$ for the smallest index $i$ such that $m_i\neq n_i$.  Then
one defines the initial (exterior) monomial with respect to $\sigma$  of any
element $f$ in the  exterior space $\wedge^d R_i$   and the   initial
subspace  of any subspace of $\wedge^d R_i$.   By construction one has that
$\ini_\sigma (V)=\langle m_1,\dots, m_d \rangle$  and $m_1>_\sigma \dots
>_\sigma m_d$ if  and only if 
$\ini_\sigma (\wedge^d V)= m_1 \wedge \dots
\wedge  m_d$.   For an element $F \in \wedge^d R_i$ we define its  
$\sigma$-support $\Support_\sigma(F)$   to be  the set of the
$\sigma$-standard exterior monomials which appear with a non-zero coefficient
in $F$. Note that any exterior monomial $n$ is equal (up to sign) to a
$\sigma$-standard exterior monomial. Note also that given an element
$F\in \wedge ^dR_i$ and two term orders $\sigma$ and $\tau$  then the
$\tau$-support of $F$ is obtained by   taking the $\tau$-standard form of
the elements in $\Support_\sigma(F)$.   One has:

\begin{lemma}
\label{lm3}
 Let  $\sigma$  be a  term order.  Let  $m=m_1
\wedge \dots \wedge  m_d$ be a $\sigma$-standard exterior monomial, and let
$q=q_1 \wedge \dots \wedge  q_d$ be an exterior monomial with 
$q_i\leq_\sigma  m_i$ for  $i=1,\dots,d$.  Let  $n=n_1 \wedge \dots \wedge  n_d$ 
be the $\sigma$-standard exterior monomial corresponding to $q$. Then
$n_i\leq_\sigma  m_i$ for  $i=1,\dots,d$.
\end{lemma} 

\begin{proof}   Since $n$ is obtained from $q$ by a sequence of
transposition exchanging $q_j$ with $q_{j+1}$  whenever
$q_j<_\sigma q_{j+1}$  it suffices to check that the property $q_i\leq_\sigma
m_i$ for all $i$ is preserved by any such a  transposition. This is clear
since 
$m_j>_\sigma m_{j+1}\geq_\sigma  q_{j+1}$ and 
$m_{j+1}\geq_\sigma  q_{j+1}>_\sigma q_j$. 
\end{proof}

\begin{lemma}
\label{lm4}
 Let  $\sigma$  be a  term order.  Let $V$ be  a
subspace of $R_i$ of dimension $d$,    and let $\ini_\sigma (\wedge^d V)=
m_1 \wedge \dots \wedge  m_d$.  For every  
$n_1\wedge\dots \wedge n_d \in \Support_\sigma(\wedge^d V)$   one has   
$m_i \geq_\sigma n_i$ for  $i=1,\dots, d$. 
\end{lemma} 

\begin{proof}  Let $f_1,\dots,f_d$ be   elements in $V$ such that  
$\ini(f_i)=m_i$. Then $\wedge^d V=f_1\wedge \dots \wedge f_d$. For   $i=1,\dots,d$
let $q_i$ be a monomial in $f_i$. It suffices to show that  the
$\sigma$-standard exterior monomial    corresponding to $q_1\wedge \dots
\wedge q_d$  satisfies the desired property. This follows from Lemma \ref{lm3} since
$q_i\leq_\sigma m_i$.   
 \end{proof}

The crucial fact is the following: 
 
\begin{lemma}
\label{lm5}
  Let $V$ be a $d$-dimensional subspace of $R_i$.
Let  
$\sigma$ and $\tau$  be   term orders. Set $W=\ini_\tau(V)$. Let $g\in
\GL_n(K)$ be a generic matrix acting as $K$-algebra graded homomorphism  on
$R$.    Then 
\[
\Support_\sigma(g(\wedge^d W))\subseteq \Support_\sigma (g(\wedge^d V)).
\]
\end{lemma}

\begin{proof} Let  $W=\langle m_1, \dots, m_d\rangle$ and let  
$f_1,\dots, f_d$ in $V$ so that $\ini_\tau(f_i)=m_i$. Set $F=f_1\wedge \dots
\wedge f_d$ and $M=m_1\wedge \dots \wedge m_d$. We have to show that
$\Support_\sigma (g(M))\subseteq  \Support_\sigma(g(F))$.   The  matrix
$g$ acts   on $R$ by, say,  $g(x_i)=\sum_j g_{ij}x_j$.  We give to
$g_{ij}$ a multidegree: $\deg(g_{ij})=e_i \in {\bf Z}^n$. In the following 
$\log(m)$ denotes the exponent of a monomial $m$.  For any monomial $m$ of
$R_i$ we have that $g(m)$ is a sum of monomials of $R_i$ whose coefficients
are polynomials of degree $\log(m)$ in the $g_{ij}$. Similarly, if 
$n=n_1\wedge \dots \wedge n_d$ is an  exterior monomial, then $g(n)$ is a sum
of $\sigma$-standard exterior monomials whose coefficients are polynomials in
the $g_{ij}$ of degree  $\log(n_1\cdots n_d)$. Now assume $n=n_1\wedge \dots
\wedge n_d$ is a $\sigma$-standard exterior monomial in the $\sigma$-support
of $g(M)$.  If
$n$ arises in the expansion of $g(Q)$ where $Q=q_1\wedge \dots \wedge q_d$
for monomials
$q_i$ in the support of $f_i$ then the coefficient of $n$ in $g(Q)$ is a
polynomial of degree $\log(q_1\dots q_d )$ in the
$g_{ij}$. If at least one of the $q_i$, say $q_j$, is $<_\tau m_i$ then
$q_1\dots q_d <_\tau m_1\dots m_d$. In particular 
$q_1\dots q_d \neq  m_1\dots m_d$. It follows that the coefficients of $n$ in
$g(M)$ and in $g(Q)$  are polynomials in the $g_{ij}$ of different degree. 
Therefore  the coefficient of $n$ in $g(F)$   is a multi-homogeneous 
polynomial in the $g_{ij}$  and one of its  homogeneous component is exactly
the coefficient of $n$ in $g(M)$.    This suffices to show that, for a
generic $g$,  the element
$n$ is in the
$\sigma$-support of $g(F)$. 
\end{proof}

\begin{corollary}
\label{cr6}
 Let $V$ be a $d$-dimensional subspace of
$R_i$.  Let  $\tau$  and $\sigma$ be  term orders.   Let
$\Gin_{\sigma}(V)=\langle m_1,\dots,m_d\rangle $ and 
$\Gin_{\sigma}(\ini_\tau(V))=\langle n_1,\dots,n_d\rangle $ with
$m_i>_\sigma m_{i+1}$ and $n_i>_\sigma n_{i+1}$ for all $i=1,\dots,d-1$. 
Then $m_i\geq_\sigma n_i$ for all $i=1,\dots,d$.
\end{corollary}

\begin{proof} Set $W=\ini_\tau(V)$,  $m=m_1\wedge \dots \wedge  m_d$ and
$n=n_1\wedge \dots \wedge  n_d$. By construction  
$\ini_\sigma(g(\wedge^d W))= n$ for a generic matrix $g$. 
By virtue of Lemma \ref{lm5}, $n \in \Support_\sigma(g(\wedge^d V))$ and by construction
$\ini_\sigma(g(\wedge^d V))=m$. It follows from Lemma \ref{lm4} that  
$n_i\leq_\sigma m_i$ for all
$i=1,\dots,d$. 
\end{proof}

We are ready to prove Theorem \ref{thm1}:

\begin{proof}[ Proof of  Theorem \ref{thm1}]  Set $H=(x_{n-p+1},x_{n-p+2},\dots,x_n)$. By
virtue of  Lemma \ref{lm2}, it is enough  to show that the Hilbert function of
$R/\Gin_{\revlex}(I)+H$ is $\leq $  that of 
$R/\Gin_{\revlex}(\ini_\tau(I))+H$.   Fix an integer $j$ and set 
\[
a=\dim\left [ R/\Gin_{\revlex}(I)+H \right ]_j, \qquad \hbox{ and } \qquad 
b=\dim\left [R/\Gin_{\revlex}(\ini_\tau(I))+H \right ]_j. 
\]
   We have to
show that $a\leq b$.  Let $V$ be the component of degree $j$ of $I$. Set
$d=\dim V$,   
$\Gin_{\revlex}(V)=\langle m_1,\dots, m_d\rangle$, 
$\Gin_{\revlex}(\ini_\tau (V))=\langle n_1,\dots, n_d\rangle$ and assume
$m_i>_{\revlex} m_{i+1}$ and  
$n_i>_{\revlex} n_{i+1}$.  For a monomial
$m$ we set
$\max(m)=\max\{ i :  x_i \hbox{ divides }  m\}$. By construction we have 
\[
b-a=|\{ k : \max(m_k)\leq n-p\}|-|\{ k : \max(n_k)\leq n-p\}|.
\]
  By
Corollary \ref{cr6} we know that 
$m_i\geq_{\revlex} n_i$ for all $i$. This implies that
$\max(m_i)\leq \max(n_i)$ for all $i$.  Hence  
$\{ k : \max(m_k)\leq n-p\} \supseteq  \{ k : \max(n_k)\leq n-p\}$ and we are
done.
\end{proof}

\begin{remark}
\label{rm7}
 One easily checks that the proof of  Theorem \ref{thm1}
works  also if one takes the initial ideal with respect to a positive weight
function $\omega$. In particular Vasconcelos' conjecture holds  in this
situation too.
\end{remark}

\begin{remark}
\label{rm8} 
With the notation of Theorem \ref{thm1}, 
 one could ask whether there is a relation between  the graded Betti
numbers  of $R/I+J$ and  those of  
$R/\ini_\tau(I)+J$.   The known properties of the initial ideal  imply that
the former are smaller than the latter for instance when 
$p\leq \depth R/\ini_\tau(I)$. But this relation does not  hold in general. 
This is because, as we know,  the Hilbert function of $R/I+J$ is $\leq $   
that of $R/\ini_\tau(I)+J$  and hence the number of generators in  low
degrees of $I+J$ tends to be larger than that $\ini_\tau(I)+J$. For
instance, taking
$I=(x^2+yz, xy, xz)$  and  $\tau$ to be the lex order, then
$\ini_\tau(I)=(x^2,xy,xz,yz^2, y^2z)$ and for a general linear form $L$ the
ideal $I+(L)$ has $3$ minimal generators in degree $2$ and
  while $\ini_\tau(I)+(L)$ has only $2$ minimal generators in degree $2$.  
\end{remark}

\begin{remark}
\label{rm9}  
Recall that the analytic spread $\ell(I)$ of an
ideal $I$ is the Krull dimension of the fiber ring 
$\oplus_{i=0}^\infty I^i/\mm I^i$. One can ask whether there is a relation
between the analytic spread $I$ and that of $\ini_\tau(I)$. There are
examples where $\ell(I)>\ell(\ini_\tau(I))$ and other where 
$\ell(I)<\ell(\ini_\tau(I))$. As for the former, take for instance the ideal
$I$  of the $2$-minors of a generic $3\times 3$ symmetric matrix and $\tau$ a
diagonal term order (i.e. the initial term of a   minor is the product of the
elements of the main diagonal). One has
$\ell(I)=6$ and
$\ell(\ini_\tau(I))=5$. On the other hand,   if $I$ is the ideal generated by
$2$ generic quadrics in
$3$ variables
 and $\tau$ is the lex order then   $\ell(I)=2$ and
$\ell(\ini_\tau(I))=3$.
\end{remark}  

\section{reduction numbers and Lex-segment ideals} 

A monomial  ideal $L$ of $R$ is  said to be a Lex-segment if whenever   $m$
is a monomial  in $L$ and   $n$ is a  monomial with $\deg(n)=\deg(m)$ and
$n>m$ with respect to the lexicographic order then  one has that $n\in L$. 
Given a homogeneous ideal
$I$ there is a unique Lex-segment ideal
$I^{\Lex}$  such that the  Hilbert function of $I^{\Lex}$ is equal to that
of $I$.    It is well-know that
$I^{\Lex}$ is ``extremal" with respect to many invariants  in the class of
the ideals with a given Hilbert function (e.g. absolute Betti numbers
 Bigatti \cite{B}, Hulett \cite{H} and Pardue \cite{P}, relative Betti numbers Iyengar and 
Pardue \cite{IP},  local cohomology Sbarra \cite{S}, etc...). Therefore it is natural to ask
whether the same holds also for the reduction number, i.e. whether  $r(R/I)\leq r(R/I^{\Lex})$
holds in general. We have:

\begin{proposition}
\label{prp10} If $K$ has  characteristic $0$,  then 
\[
r(R/I)\leq
r(R/I^{\Lex})
\]
 holds for every homogeneous ideal of  $I$ of
$K[x_1,\dots,x_n]$. 
\end{proposition}

\begin{proof}  Let $I$ be a homogeneous ideal of $R=K[x_1,\dots, x_n]$ and
set $d=\dim R/I$ and 
$J=\Gin_{\revlex}(I)$.   It is know that
$J$ is Borel-fixed, that is fixed under the action of the group of the
upper-triangular matrix. In characteristic $0$ this is equivalent to say
that $J$ is  strongly stable, that is, if
$x_im$ is a monomial in $J$ and $j<i$ then $x_jm$ is in $J$ as well.  Form
this and from  Lemma \ref{lm2}, it follows immediately that if $\chara K=0$, then  
$r(R/I)$ is equal to the least integer $k$ such that $x_{n-d}^{k+1}$ is in
$J$ (this  has been observed also in \cite{T}).  Then the desired inequality is a
consequence of the following fact: 

\noindent{ Claim}: Let  $V$ and $L$ be  sets  of monomials of degree $k$ with the
same cardinality such that $V$ is  strongly stable and $L$ is a Lex-segment.
If $x_i^k\in L$ for some
$i$,     then $x_i^k\in V$. 

To prove the claim one  observes that since $L$ contains $x_i^k$ and it is a
Lex-segment, then $L$ contains also the set, say $A$, of all the monomials
of degree $k$ which are divisible by some $x_j$ with $j<i$. Therefore
$|L|\geq |A|+1$. Since $|V|=|L|$ it follows that $V$   contains a monomial
$m$ which is not in $A$.  In other words, $V$  contains a monomial supported
only on the variables $x_i,x_{i+1},\dots,x_n$. Since $V$ is stable, then $V$
contains also $x_i^k$.\end{proof}

We  believe that the inequality  of Proposition \ref{prp10} is   true  also if the
characteristic of the base field is finite.  Pardue developed in \cite{P}   a 
characteristic free  strategy for proving that the Lex-segment ideal is
extremal with respect to a certain invariant, say 
$\alpha(R/I)$.  Roughly speaking, it says that if $\alpha$ does not decrease
by taking initial ideals and also does not decrease by performing a certain
deformation process, called  polarization,  then one has 
$\alpha(R/I)\leq
\alpha(R/I^{\Lex})$ for all the ideals $I$.  For the  definition of  
polarization of a monomial ideal we refer the reader to the paper of Pardue
\cite{P}. Unfortunately one cannot use Pardue's  argument to prove the above
inequality.    This is because the reduction number can decrease under
polarizations. For example, let
$R=K[x_1,\dots, x_4]$ and  
\[
I=(x_4^2, x_1x_3^3, x_3^3x_4, x_2^3x_4, x_2x_3^3, x_2^3x_3, x_1^2x_3^2,
x_1^4, x_1x_2^2x_4, x_2^4)
\]
 and $J$ its polarization; one can check that
$r(R/I)=4$ and    $r(R/J)=3$. In this case $r(R/I^{\Lex})=5$. 

\medskip

Let us also note that the above ideal  can be used to construct an example
of a standard graded algebra $A$ and a non-zero divisor $z$ of degree $1$
such that $r(A)<r(A/zA)$. To this end it suffices to take $S=K[x_1,\dots,
x_5]$,  and  
\[
I_1=(x_4x_5, x_1x_3^3, x_3^3x_4, x_2^3x_4, x_2x_3^3, x_2^3x_3, x_1^2x_3^2,
x_1^4, x_1x_2^2x_4, x_2^4).
\]
In other words, $I_1$ is the polarization of
the ideal $I$ above with respect to the variable
$x_4$. Set
$A=S/I_1$ and
$z=x_4-x_5$.  Then
$z$ is a non-zero divisor of
$A$ and  $r(A)=3$ and  $r(A/zA)=4$.

\medskip 

\noindent{\bf Thanks:}\   {We would like to thank Marilina Rossi   for useful
discussions concerning the material of this paper.  The explicit examples  
that we have  presented  in the paper have been detected by  using  the
computer algebra system CoCoA  \cite{CNR}.
}

\medskip 

\bibliographystyle{amsplain}

\end{document}